\documentclass[leqno]{amsart}
\usepackage[cp1250]{inputenc}
\usepackage[T1]{fontenc}
\usepackage{amsthm}
\usepackage[leqno]{amsmath}
\usepackage{amssymb}
\usepackage{amsfonts}
\usepackage{enumerate}
\usepackage{latexsym}

\makeatletter
\@namedef{subjclassname@2010}{%
  \textup{2010} Mathematics Subject Classification}
\makeatother

\newcommand{\kk}{\mathbf{\mathbb{K}}}

\newcommand{\rr}{\mathbf{\mathbb{R}}}
\newcommand{\cc}{\mathbf{\mathbb{C}}}
\newcommand{\nn}{\mathbf{\mathbb{N}}}

\newcommand{\ord}{\operatorname{ord}}

\newcommand{\dist}{\operatorname{dist}}

\newcommand{\wykl}{{\mathcal{L}}}
\newcommand{\LL}{\mathbf{L}}

\newcommand{\er}{\mathbb{R}}

\newcommand{\en}{\mathbb{N}}
\newcommand{\ez}{\mathbb{Z}}

\newcommand{\parallell}{||}
\newcommand{\parallelp}{||}
\newcommand{\ec}{\mathbb{C}}

\newcommand{\wym}{\,{\rm dim}\,}

\newcommand{\graph}{\,{\rm graph}\,}

\newcommand{\wl}{{\mathcal{L}}_{\infty }}

\newcommand{\Ld}{\Delta}

\def\Int{\operatorname{Int}}

\def\wl{{\mathcal{L}}_{\infty}}

\newtheorem{theorem}{Theorem}[section]
\newtheorem{lemma}[theorem]{Lemma}

\newtheorem{cor}[theorem]{Corollary}
\newtheorem{proposition}[theorem]{Proposition}

\theoremstyle{remark}

\numberwithin{equation}{section}

\title{Metric properties of semialgebraic mappings}

\author[K. Kurdyka]{Krzysztof Kurdyka} 
\address{Laboratoire de Mathematiques (LAMA) Universi\'e de Savoie,  UMR-5127 de CNRS\newline
73-376 Le Bourget-du-Lac cedex FRANCE}
\email{Krzysztof.Kurdyka@univ-savoie.fr}

\author[S. Spodzieja]{Stanis{\l}aw Spodzieja}
\address{Faculty of Mathematics and Computer Science, University of \L \'od\'z\newline
S. Banacha 22, 90-238 \L \'od\'z, Poland}
\email{spodziej@math.uni.lodz.pl}

\author[A. Szlachci\'nska]{Anna Szlachci\'nska}
\address{Faculty of Mathematics and Computer Science, University of \L \'od\'z\newline
S. Banacha 22, 90-238 \L \'od\'z, Poland}
\email{anna$\_$loch@wp.pl}

\date{\today}

\begin{document}

\begin{abstract}We give an effective estimation from above for the local {\L}ojasiewicz exponent for separation of semialgebraic sets and for a semialgebraic mapping  on a closed semialgebraic set. We also give an effective estimation from below of the {\L}ojasiewicz exponent in the global separation for semialgebraic sets and estimation of the {\L}ojasiewicz exponent at infinity of a semialgebraic mapping similar to the Jelonek result \cite{J2} in the complex case. Moreover, we prove that both local and  global {\L}ojasiewicz exponent of an overdetermined semialgebraic mapping $F:X\to \rr^m$  on a closed semialgebraic set $X\subset \rr^n$ (i.e. $m>\dim X$) are equal to the {\L}ojasiewicz exponent of the composition $L\circ F:X\to\rr^n$ for the generic linear mapping $L:\rr^m\to \rr^k$, where $k=\dim X$.
\end{abstract}

\subjclass[2010]{14P20, 14P10, 32C07.} 
\keywords{{\L}ojasiewicz exponent, semialgebraic set, semialgebraic mapping, polynomial mapping, overdetermined mapping.}

\maketitle

\section*{Introduction}

{\L}ojasiewicz inequalities are an important and useful tool in differential equations, singularity theory and optimization (see for instance 
\cite{KMP} in the local case and 
\cite{RS3} \cite{RS4} 
at infinity).  In these considerations, an estimations of the local and global {\L}ojasiewicz exponents (see for instance 
\cite{K2},   \cite{LejeuneTeissier}, \cite{RS2}, \cite{RS3}, \cite{T}  in the local case and  
\cite{Cy}, \cite{CKT}, 
\cite{JKS}, \cite{K1}, \cite{RS1} at infinity) play a central role. 
In the complex case, an essential estimations of the {\L}ojasiewicz exponent at infinity of a polynomial mapping $F=(f_1,\ldots,f_m):\cc^N\to \cc^m$ on an algebraic set $V\subset \cc^N$ (see Section~\ref{Lojasiewiczexpatinfty}) denoted by $\wykl_\infty^\cc(F|V)$ or $\wykl_\infty^\cc(F)$ for $V=\cc^n$, was obtained by J.~Ch\k{a}dzy\'nski \cite{CK3}, J.~Koll\'ar \cite{K1}, E. Cygan, T. Krasi\'nski, P. Tworzewski \cite{CKT}, Z.~Jelonek \cite{J1}, \cite{J2} and  E. Cygan \cite{Cy}. 
More precisely, let $\deg f_j=d_j$, $j=1,\ldots,m$, $d_1\ge \ldots\ge d_m>0$ and let
$$
B(d_1,\ldots,d_m;k)=\begin{cases}d_1\cdots d_m&\hbox{for}\;\;m\leqslant k,\\
d_1\cdots d_{k-1}d_m&\hbox{for}\;\;m>k.\end{cases}
$$
J.~Ch\k{a}dzy\'nski proved that 
\begin{equation}\label{Chadzynski}\tag{Ch}
\wykl^\cc_\infty(F)\ge d_2 -d_1d_2+\sum_{b\in F^{-1}(0)}\mu_b(F),
\end{equation} where $\mu_b(F)$ is the multiplicity of $F$ at $b$, provided $N=m=2$ and $\#F^{-1}(0)<\infty$.
For arbitrary $m\ge N$, under the assumption $\#F^{-1}(0)<\infty$, J. Koll\'ar proved that
\begin{equation}\label{Kollar}\tag{K}
\wl^\cc(F)\ge d_m-B(d_1,\ldots,d_m;N),
\end{equation}
and E. Cygan, T. Krasi\'nski, P. Tworzewski gave an estimate
\begin{equation}\label{CyganKrasinskiTworzewski}\tag{CKT}
\wl^\cc(F)\geq d_m-B(d_1, \ldots, d_m;N)+\sum_{b \in
F^{-1}(0)}\mu_b(F),
\end{equation}
where $\mu_b(F)$ is the intersection multiplicity in the sense of R. Achilles, P.~Tworzewski and  T. Winiarski of $\graph F$ and $\cc^n\times \{0\}$ at the point $(b,0)$ (see \cite{ATW}). Z.~Jelonek obtained the following estimation for a complex $k$-dimmensional algebraic variety $V\subset \cc^N$ of degree $D$:
\begin{equation}\label{Jelonek}\tag{J}
\wl^\cc(F|V)\geq d_m-D\cdot B(d_1, \ldots,d_m;k)+\nu,
\end{equation}
where $\nu= \#(F^{-1}(0)\cap V)<+\infty$. E. Cygan gave the following global 
inequality 
\begin{equation}\label{ECygan}\tag{C$_1$}
|F(x)|\geq C \left(\frac{\dist(x,F^{-1}(0))}{1+|x|^2}\right)^{B(d_1, \ldots,
d_m;N)}\quad {\rm for} \quad x \in \cc^N
\end{equation}
for some positive constant $C$. Moreover, E. Cygan proved that for a complex algebraic sets $X,Y\subset \cc^N$ there exists a positive constant $C$ such that
\begin{equation}\label{ECygan2}\tag{C$_2$}
\dist (x,X)+\dist (x,Y)\ge C\left(\frac{\dist(x,X\cap Y)}{1+|x|^2}\right)^{\deg X\cdot\deg Y}\quad\hbox{for } x\in\cc^N.
\end{equation}
A similar result to \eqref{ECygan2} was obtained by S. Ji, J. Koll\'ar, B. Shiffman \cite{JKS}.

For the real algebraic sets we have the following global {\L}ojasiewicz inequality 
(see \cite
{KS}): 
 if $X, Y\subset \rr^N$ are algebraic sets defined by a systems of polynomial equations of degrees at most $d$, then for some positive constant $C$,
\begin{equation}\label{eqKSmain}\tag{KS$_1$}
\dist(x,X)+\dist(x,Y)\ge C\left(\frac{\dist(x,X\cap Y)}{1+|x|^2}\right)^{d(6d-3)^{N-1}}\;\;\hbox{for }x\in \rr^N.
\end{equation}
In particular, we have the following global {\L}ojasiewicz inequality 
 (see \cite
{KS}). 
Let $F=(f_1,\ldots,f_m):\er^N\to\er^m$ be a polynomial mapping of degree $d$. Then for some positive  constant $C$,
\begin{equation}\label{eq8KS}\tag{KS$_2$}
|F(x)|\ge C\left(\frac{\dist(x,F^{-1}(0))}{1+|x|^2}\right)^{{d(6d-3)^{N-1}}} \quad\hbox{for}\quad x\in\er^N.
\end{equation}
If, additionally, the set $F^{-1}(0)$ is compact, then  
\begin{equation}\label{eq83}\tag{KS$_3$}
\wykl^\rr_\infty(F)\ge {-d(6d-3)^{N-1}}.
\end{equation}

The purpose of this article is to generalize a  results similar to \eqref{Jelonek}, \eqref{ECygan}, \eqref{ECygan2} and \eqref{eqKSmain}, \eqref{eq8KS}, \eqref{eq83} in the case of algebraic sets and regular mappings to the case of semialgebraic sets and mappings. More precisely, we give an effective estimation from above for the local {\L}ojasiewicz exponent for separation of semialgebraic sets (see Theorem \ref{thmsemialglocsets} in Section \ref{Lojasiewiczexpatapoint})  and for a semialgebraic mapping  on a closed semialgebraic set (see Corollaries \ref{thmsemialglocal} and \ref{estlocalonalgebraicser} in Section \ref{Lojasiewiczexpatapoint}). We also give an effective estimation from below of the {\L}ojasiewicz exponent in the global separation for semialgebraic sets (see Theorem \ref{thmsemialglocsetsinfty} in Section \ref{Lojasiewiczexpatinfty}) and estimation of the {\L}ojasiewicz exponent at infinity of a semialgebraic mapping similar to the Jelonek result \cite{J2} in the complex case (see Corollaries \ref{thmsemialg}, \ref{thmsemialg} in Section \ref{Lojasiewiczexpatinfty}). The above estimations are effective in terms of degrees of 
polynomials describing semialgebraic sets and mappings. Moreover, we prove that both: local and  global {\L}ojasiewicz exponent of an overdetermined analytic and semialgebraic mapping $F:X\to \rr^m$  on a closed semialgebraic set $X\subset \rr^N$ (i.e. $m>\dim X$) are equal to the {\L}ojasiewicz exponent of the composition $L\circ F:X\to\rr^n$ for the generic linear mapping $L:\rr^m\to \rr^k$, where $k=\dim X$ (see Theorems \ref{Theorem1}, \ref{Theorem2},  and \ref{Theorem3}). 



A mapping $F:\kk^N\to\kk^m$ is called \emph{overdetermined} if $m>N$. These results are obtained by extending the mappings to overdetermined ones and reduction of calculations for the both local and global {\L}ojasiewicz exponents of overdetermined polynomial mappings to the case when $m=N$ (cf \cite{Sp2}, \cite{Sp3}). The crucial roles in the proofs are played by inequalities \eqref{eqKSmain} and \eqref{eq8KS}.

\section{The {\L}ojasiewicz exponent at a point}\label{Lojasiewiczexpatapoint}

We will give an estimate from above of the {\L}ojasiewicz exponent for the regular se\-pa\-ra\-tion of closed semialgebraic sets and for a continuous semialgebraic mapping on a closed semialgebraic set. 
Let us start from some notations. Let $X\subset \rr^N$ be a closed semialgebraic set. It is known that $X$ has the decomposition \begin{equation}\label{decompositionsemialg} 
X=X_1\cup\ldots\cup X_k
\end{equation}
into the sum of closed basic semialgebraic sets  
\begin{equation}\label{semialgebrdecomposition}
X_i=\{x\in\rr^N:g_{i,1}(x)\ge 0,\ldots,g_{i,r_i}(x)\ge 0,\; h_{i,1}(x)=\cdots=h_{i,l_i}(x)=0\}, 
\end{equation}
$ i=1,\ldots k$ (see \cite{BCR}), where $g_{i,1},\ldots,g_{i,r_i},  h_{i,1},\ldots ,h_{i,l_i}\in\rr[x_1,\ldots,x_N]$. 
Assume that $r_i$ is the smallest possible number of inequalities  $g_{i,j}(x)\ge 0$ in the definition of $X_i$, for $i = 1, ...,k$. 
Denote by $r(X)$ the minimum of $\max\{r_1,\ldots,r_k\}$ 
for any decomposition \eqref{decompositionsemialg} into sum of sets of form \eqref{semialgebrdecomposition}. As  L. Br\"ocker \cite{Brocker0} (cf., \cite{Brocker1}, \cite{Scheiderer}) showed, 
\begin{equation*} 
r(X)\le \frac{1}{2}N(N+1).
\end{equation*}
 Denote by $\kappa(X)$ the mimimum of numbers 
$$
\max\{\deg g_{1,1},\ldots,\deg g_{k,r_k}, \deg h_{1,1},\ldots ,\deg h_{k,l_k}\}
$$ 
for any decomposition \eqref{decompositionsemialg} of $X$ into the sum of sets of form \eqref{semialgebrdecomposition}, provided $r_i\le r(X)$. Obviously $r(X)=0$ if and only if $X$ is an algebraic set. The numbers $r(X)$ and $\kappa(X)$ characterize the so called complexity of semialgebraic set $X$. For more information about the complexity see for example \cite{BPR1}, \cite{BPR2}, \cite{BPR}, \cite{BCR}, \cite{Roy}.

\begin{theorem}\label{thmsemialglocsets}
Let $X,Y\subset \rr^N$ be a closed semialgebraic sets, and let $0\in X\cap Y$. Put $r=r(X)+r(Y)$ and $d=\max\{\kappa(X),\kappa(Y)\}$.  Then there exist a neighbourhood $U\subset \rr^N$ of $0$ and some positive constant $C$, such that
\begin{equation}\label{eqsemialgsets1local}
\dist (x,X)+\dist (x,Y)\ge C\dist (x,X\cap Y)^{d(6d-3)^{N+r-1}}\;\;\hbox{for }x\in {U}.
\end{equation}
If, additionally, $0$ is an isolated point of $X\cap Y$, then for some neighbourhood $U\subset \rr^N$ of $0$ and some positive constant $C$,
\begin{equation}\label{eqsemialgsets2local}
\dist (x,X)+\dist (x,Y)\ge C|x|^{\frac{(2d-1)^{N+r}+1}{2}}\;\;\hbox{for }x\in {U}.
\end{equation}
\end{theorem}

The proof of the above theorem will be carried out in Section~\ref{sectsemialgebrlocal}. The key role in the proof will be played the following inequality 
\cite[Corollary 8]{KS}: let  $X=(g_1,\ldots,g_k)^{-1}(0)$ and $Y=(h_1,\ldots,h_l)^{-1}(0) \subsetneq \er^N$, where  $g_i,h_j\in \rr[x_1,\ldots,x_N]$ are polynomials  of degree not greater than $d$. Let $a\in \er^N$.
Then there exists a positive constant $C$ such that
\begin{equation}\label{eq7}\tag{KS$_4$}
\dist (x,X)+\dist (x,Y)\ge C\dist (x,X\cap Y)^{d(6d-3)^{N-1}}
\end{equation}
in a neighbourhood of  $a$. If, additionally, $a$ is an isolated point of $X\cap Y$, then 
\begin{equation}\label{eq71}\tag{KS$_5$}
\dist (x,X)+\dist (x,Y)\ge C |x-a|^{\frac{(2d-1)^{N}+1}{2}}
\end{equation}
in a neighbourhood of  $a$ for some positive $C\in\rr$.

\medskip
Let $X\subset \kk^N$ be a closed subanalytic set. If $\kk=\cc$ we consider $X$ as a subset of $\rr^{2N}$. We will assume that the origin $0\in \kk^N$ belongs to $X$ and it is an accumulation point of $X$. We denote by $F:(X,0)\to (\kk^m,0)$ a mapping of a neighbourhood $U\subset X$ of the point $0\in\kk^N$ into $\kk^m$ such that $F(0)=0$, where the topology of $X$ is induced from $\kk^N$. 

Let $F:(X,0)\to (\kk^m,0)$ be a continuous subanalytic mapping, i.e. the graph of $F$ is a closed subanalytic subset of $(X\cap U)\times \kk^m$ for some neighbourhood $U\subset \kk^N$ of the origin. If $\kk=\cc$, we consider $\kk^N$ as $\rr^{2N}$ and $\kk^m$ as $\rr^{2m}$. 
Then there are positive constants $C,\eta,\varepsilon$ such that the following \emph{{\L}ojasiewicz inequality} holds:
\begin{equation}\label{Lojineq0}
|F(x)|\ge C\dist(x,F^{-1}(0)\cap X)^\eta\quad \hbox{if}\quad x\in X,\quad |x-a|<\varepsilon,
\end{equation}
where  $|\,\cdot\,|$ is the Euclidean norm in $\kk^n$, respectively in $\kk^N$, and $\dist (x,V)$ is the distance of  $x\in\kk^N$ to the set $V\subset \kk^N$ ($\dist (x,V)=1$ if $V=\emptyset$). The smallest exponent $\eta$ in \eqref{Lojineq0} is called the \emph{{\L}ojasiewicz exponent} of $F$ \emph{on the set} $X$ \emph{at} $0$ and is denoted by $\wykl^\kk_0(F|X)$. If $X$ contains a neighbourhood $U\subset \kk^N$ of $0$ we will call it the {\it {\L}ojasiewicz exponent} of $F$ \emph{at} $0$  and denote by $\wykl^\kk_0(F)$. It is known that $\wykl^\kk_0(F|X)$ is a rational number and \eqref{Lojineq0} holds with any $\eta\ge\wykl^\kk_0(F|S)$ and some positive constants $C,\varepsilon$, provided $0$ is an accumulation point of $X\setminus F^{-1}(0)$ (see \cite{BR}, \cite{Sp4}). If $0$ is an isolated piont of $X\setminus F^{-1}(0)$, we have $\wykl^{\kk}_0(F|X)=0$. 

\medskip
From Theorem \ref{thmsemialglocsets} there follows

\begin{cor}\label{thmsemialglocal}
Let $F:X\to\rr^m$ be a continuous semialgebraic mapping, where $X\subset \rr^N$ is a closed semialgebraic set, let $0\in X$ and $F(0)=0$. Put $r=r(X)+r(\graph F)$ and $d=\max\{\kappa(X),\kappa(\graph F)\}$. Then 
\begin{equation}\label{eqsemialg1local}
\wykl^\rr_0(F|X)\le d(6d-3)^{N+r-1}.
\end{equation}
If, additionally, $0$ is an isolated zero of $F$, then 
\begin{equation}\label{eqsemialg1localx}
\wykl^\rr_0(F|X)\le \frac{(2d-1)^{N+r}+1}{2}.
\end{equation}
\end{cor}

\medskip
For a real polynomial mapping $F:\rr^N\to\rr^m$ a similar result as above was obtained in 
 \cite[Corollary 6]{KS}. Namely, if $d=\deg F$, then 
\begin{equation}\label{KurdykaSpodz}\tag{KS$_6$}
\wykl_0^\rr(F)\le d(6d-3)^{N-1}.
\end{equation}
In the case of regular mapping, i.e. the restriction of polynomial mapping to algebraic set, from Corollary \ref{thmsemialglocal}, 
Theorem \ref{Theorem2}  (see Section \ref{Genericsection}) and \cite{Cy2} (see also \cite{CKT}, \cite{Cy}) we obtain an estimation of its local {\L}ojasiewicz exponent, also  for regular mapping with nonisolated zero-set (cf \cite{P2}, \cite{Sp3} for mapping with isolated zeroes).

\begin{cor}\label{estlocalonalgebraicser}
Let $F:(\kk^N,0)\to(\kk^m,0)$, $m\ge N$,  be a polynomial mapping, let $X\subset \kk^N$ be an algebraic set defined by a system of  equations $g_1(x)= \ldots =g_r(x)=0$, where $g_1,\ldots, g_r\in\kk[x_1,\ldots,x_N]$, and let $d=\max\{\deg F,\deg g_1,\ldots,\deg g_r\}$. Assume that $d>0$ and $0\in X$.

{\rm(a)} If $\kk=\rr$, then $\wykl_0^\rr(F|X)\le d(6d-3)^{N-1}$.

{\rm (b)} If $\kk=\cc$, then $\wykl_0^\cc(F|X)\le d^N$.
\end{cor}

Indeed, the assertion (a) immediately follows from Corollary \ref{thmsemialglocal}. We will prove the assertion (b). Let $G=(F,g_1,\ldots,g_r):\cc^N\to \cc^{m+r}$. By Theorem \ref{Theorem2}, for the generic $L=(L_1,\ldots,L_m)\in
\Ld^\cc(m+r,m)$ we have $\wykl_0^\cc(G|X)=\wykl_0^\cc(L\circ G|X)$. Moreover, $\deg L_j\circ G\le d$ for $j=1,\ldots,m$. E.~Cygan in \cite{Cy2}  proved that for analytic sets $Z,Y\subset \cc^{N+m}$ the intersection index at $0$ of $Z$ and $Y$ is a separation exponent of $Z$ and $Y$ at the point $0\in Z\cap Y$. 
It is known that for $Z=\cc^N\times \{0\}$ and $Y=\graph L\circ G$, the index does not exceed $d^N$ (see \cite {Tw}, \cite{CKT}), so, $\wykl_0^\cc(L\circ G)\le d^N$. Since $G^{-1}(0)=F^{-1}(0)\cap X$ and by the definition of $\Ld^\cc(m+r,m)$, we have $G(x)=(F(x),0)$ for $x\in X$, then $\wykl_0^\cc(F|X)\le d^N$. This gives the assertion (b).

\section{The {\L}ojasiewicz exponent at infinity}\label{Lojasiewiczexpatinfty}

The second aim of this article is to obtain a similar results as in the previous section but for the {\L}ojasiewicz exponent at infinity. 

By the \emph{{\L}ojasiewicz exponent at infinity of a mapping $F:X\to\kk^m$} we mean the supremum of the exponents $\nu$ in the following \emph{{\L}ojasiewicz inequality}:
\begin{equation}\label{eq82}
|F(x)|\ge C|x|^\nu\quad\hbox{for}\quad x\in X,\quad |x|\ge R
\end{equation}
for some positive constants $C$, $R$; we denote it by $\wykl^\kk_\infty(F|X)$. If $X=\kk^N$ we call the exponent $\wykl^\kk_\infty(F|X)$ the \emph{{\L}ojasiewicz exponent at infinity} of $F$ and denote by $\wykl^\kk_\infty(F)$.

\medskip
By using \eqref{eq8KS} we obtain a global {\L}ojasiewicz inequality for regular mappings. 

\begin{cor}\label{estimationonanalgebraic} 
Let $X\subset \rr^N$ be an algebraic set defined by a system of polynomial equations $g_1(x)=\cdots=g_r(x)=0$, where $g_1,\ldots, g_r\in \rr[x_1,\ldots,x_N]$. Let $F:\rr^N\to\rr^m$ be a polynomial mapping and let $d=\max\{\deg F, \deg g_1,\ldots,\deg g_r\}$. 
Then for some positive constant $C$,
$$
|F(x)|\ge C\left(\frac{\dist(x,F^{-1}(0)\cap X)}{1+|x|^2}\right)^{{d(6d-3)^{N-1}}} \quad\hbox{for}\quad x\in X.
$$
If, additionally, $X$ is an unbounded set and $F^{-1}(0)\cap X$ is a compact set, then
$$
\wykl^\rr_\infty (F|X)\ge - d(6d-3)^{N-1}.
$$
\end{cor}

Indeed, let $G=(g_1,\ldots, g_r):\rr^N\to\rr^r$, and  
let $H:\rr^N\to\rr^{m+r}$ be a polynomial mapping defined by $H(x)=(F(x),G(x))$ for $x\in \rr^N$. Then $H^{-1}(0)=F^{-1}(0)\cap X$ so, by \eqref{eq8KS} we deduce the first part of the assertion. If $F^{-1}(0)\cap X$ is a compact set, then $H^{-1}(0)$ is a compact set, too. So, the second part of the assertion follows immediately from the first one (cf \eqref{eq83}). 

In the above proof we cannot apply \eqref{Chadzynski}, \eqref{Kollar}, \eqref{CyganKrasinskiTworzewski}, \eqref{Jelonek} and \eqref{ECygan}, 
 because the complexification of a real regular mapping with compact zero-set can have an unbounded zero-set.

In Section~\ref{sectsemialgebr} we will prove the following global {\L}ojasiewicz inequality for semialgebraic sets.

\begin{theorem}\label{thmsemialglocsetsinfty}
Let $X,Y\subset \rr^N$ be a closed semialgebraic sets. Put $r=r(X)+r(Y)$ and $d=\max\{\kappa(X),\kappa(Y)\}$. 
 Then there exists a positrive constant $C\in\rr$ such that
\begin{equation}\label{eqsemialgsets1infty}
\dist (x,X)+\dist (x,Y)\ge C\left(\frac{\dist (x,X\cap Y)}{1+|x|^d}\right)^{d(6d-3)^{N+r-1}}\quad\hbox{for }x\in\rr^N.
\end{equation}
\end{theorem}

The crucial role in the proof of the above theorem is played by \eqref{eqKSmain}.

Since for semialgebraic mapping $F:X\to \rr^m$, $X\subset \rr^N$ the complement of the set $(X\times \{0\})\cup \graph F\subset \rr^{N+m}$ is a dense subset of $\rr^{N+m}$, from Theorem \ref{thmsemialglocsetsinfty} immediately implies

\begin{cor}\label{thmsemialg}
Let $F:X\to\rr^m$ be a continuous semialgebraic mapping, where $X\subset \rr^N$ is a closed semialgebraic set. If $d=\max\{2,\kappa(X),\kappa(Y)\}$ 
and $r=r(X)+r(Y)$, where $Y=\graph F$, then there exists a positive constant $C\in\rr$ such that 
\begin{equation}\label{eqsemialg1KS2}
|F(x)|\ge C\left(\frac{\dist(x,F^{-1}(0)\cap X)}{1+|x|^{d}}\right)^{{d(6d-3)^{N+r-1}}} \quad\hbox{for}\quad x\in X.
\end{equation}

In particular, if $X$ is an unbounded set and $F^{-1}(0)\cap X$ is a compact set, then
\begin{equation}\label{eqsemialg12}
\wykl^\rr_\infty(F|X)\ge (1-d)d(6d-3)^{N+r-1}.
\end{equation}
\end{cor}

For a polynomial mapping $F:X\to \rr^m$ we have $r(\graph F)=r(X)$ and $\kappa(\graph F)=\max\{\deg F,\kappa(X)\}$, so we have

\begin{cor}\label{thmsemialg2}
Let $F:X\to\rr^m$ be a polynomial mapping, where $X\subset \rr^N$ is a closed semialgebraic set. If $D=\max\{2,\kappa(X)\}$, and $d=\max \{\deg F,D\}$, and $r=2r(X)$, then 
\begin{equation}\label{eqsemialg1KS22}
|F(x)|\ge C\left(\frac{\dist(x,F^{-1}(0)\cap X)}{1+|x|^{D}}\right)^{{d(6d-3)^{N+r-1}}} \quad\hbox{for}\quad x\in X.
\end{equation}

In particular, if $X$ is an unbounded set and $F^{-1}(0)\cap X$ is a compact set, then
\begin{equation}\label{eqsemialg122}
\wykl^\rr_\infty(F|X)\ge -\frac{D}{2}d(6d-3)^{N+r-1}.
\end{equation}
\end{cor}

The above corollary is not a direct consequence of Corollary \ref{thmsemialg}, so we will prove it in Section \ref{sectsemialgebrlocal}.

\section{Composition of semialgebraic mapping\\ with generic linear mapping}\label{Genericsection}

Let $\kk=\rr$ or $\kk=\cc$. By the dimension $\dim_{\kk,0} X$ at $0$ of a set $X\subset \kk^N$ we mean the infimum of the dimensions over $\kk$ at $0$ of local analytic sets $0\in V\subset \kk^N$ such that $X\cap U\subset V$ for some neighbourhood $U\subset \kk^n$ of $0$. 

By the dimension $\dim_\rr X$ of a set $X\subset \rr^N$ we mean the infimum of dimensions of local analytic sets $V\subset \rr^N$ such that $X\subset V$. In particular, if $X$ is a semialgebraic set, $\dim_\rr X$ is the infimum of dimensions of algebraic sets $V\subset \rr^N$ such that $X\subset V$ and there exists a ball $B\subset \rr^N$ centered at  $0$ such that $\dim _\rr (X\cap B)=\dim_{\rr,0} X$.

We will write "for the generic $x\in A$" instead of "there exists an
algebraic set $V$ such that $A\setminus V$ is a dense subset of $A$ and for $x\in 
A\setminus V$".  

By $\LL^\kk(m,k)$ we shall denote the set of all linear mappings
$\kk^m\to\kk^k$ (we identify $\kk^0$ with $\{0\}$). Let $m\ge k$. By $\Ld^\kk(m,k)$ we denote the set of all linear mappings
$L \in \LL^\kk(m,k)$ of the form $L =(L_1,...,L_k)$, 
$$
L_i(y_1,...,y_m)=y_i+\sum_{j=k+1}^m\alpha_{i,j} y_j,\qquad i=1,...,k, 
$$
where $\alpha_{i,j}\in \kk$.

In Section \ref{Proofssection} we will prove 
(cf \cite[Theorem 2.1]{Sp3} and \cite[Theorem 1]{SS})

\pagebreak
\begin{theorem}\label{Theorem1} Let $F=(f_1,\ldots,f_m):(X,0)\to(\rr^m,0)$ be an analytic mapping with isolated zero at the origin, where $X\subset \rr^N$ is a closed semialgebraic set and $0\in X$. Let $\dim_{\rr,0} X=n$, and let $n\le k\le m$. Then for any $L\in \LL^\rr(m,k)$ such
that the origin is an isolated zero of $L\circ F|X$, we have 
\begin{equation}\label{eq1realisolated}
\wykl^\rr_0(F|X)\le\wykl^\rr_0(L\circ F|X).
\end{equation}
Moreover, for the generic $L\in \LL^\rr(m,k)$ the origin is an isolated zero of $L\circ F|X$ and 
\begin{equation}\label{eq2realisolated}
\wykl_0^\rr(F|X)=\wykl_0^\rr(L\circ F|X).
\end{equation}
In particular, for the generic $L\in
\Ld^\rr(m,k)$ the origin is an isolated zero of $L\circ F|X$ and \eqref{eq2realisolated} holds. 
\end{theorem}

The above theorem gives a method for reduction of the problem of calculating the {\L}ojasiewicz exponent of overdetermined mappings to the case where the dimensions of domain and counterdomain of mappings are equal. It is not clear to the authors whether the above statement is true if the set $ X $ is subanalytic instead of semialgebraic.

If $F:X\to \kk^m$ is a semialgebraic mapping then without any assumptions on the set of zeroes of $F$ we will prove in Section \ref{Proofssection2} the following 

\begin{theorem}\label{Theorem2} Let $F:(X,0)\to(\kk^m,0)$ be a  continuous semialgebraic mapping, $X\subset \kk^N$ be a closed semialgebraic set of dimension $\dim_{\rr,0} X=n$, and let $n\le k\le m$. Then for any $L\in \LL^\kk(m,k)$ such that
\begin{equation}\label{condition1}
F^{-1}(0)\cap U_L=(L\circ F)^{-1}(0)\cap U_L\hbox{ for a neighbourhood $U_L\subset X$ of\; $0$}
\end{equation}
 we have 
\begin{equation}\label{eq1lokalnonisolated}
\wykl^\kk_0(F|X)\le\wykl^\kk_0(L\circ F|X).
\end{equation}
Moreover, for the generic $L\in \LL^\kk(m,k)$ the condition \eqref{condition1} holds and
\begin{equation}\label{eq2lokalnonisolated}
\wykl_0^\kk(F|X)=\wykl_0^\kk(L\circ F|X).
\end{equation}
In particular, for the generic $L\in \Ld^\kk(m,k)$ hold \eqref{condition1} and \eqref{eq2lokalnonisolated}. 
\end{theorem}

In Section \ref{ProofofTheorem} we will prove the following version of Theorem \ref{Theorem1} for the {\L}ojasiewicz exponent at infinity (cf \cite[Theorem 2.1]{Sp2}, \cite[Theorem 3]{SS}).

\begin{theorem}\label{Theorem3} Let $F=(f_1,\ldots,f_m):X\to\rr^m$ be a continuous semialgebraic mapping having a compact set of zeros, where $X\subset \rr^N$ is a closed semialgebraic set, $\dim X=n$, and let $n\le k\le m$. Then for any $L\in \LL^\rr(m,k)$ such that\linebreak $(L\circ F)^{-1}(0)\cap X$ is compact, we have 
\begin{equation}\label{eq1global}
\wykl^\rr_\infty(F|X)\ge\wykl^\rr_\infty(L\circ F|X).
\end{equation}
Moreover, for the generic $L\in \LL^\kk(m,k)$ the set $(L\circ F)^{-1}(0)$ is compact and
\begin{equation}\label{eq2global}
\wykl_\infty^\rr(F|X)=\wykl_\infty^\rr(L\circ F|X).
\end{equation}
In particular, \eqref{eq2global} holds for the generic $L=(L_1,...,L_k)\in
\Ld^\rr(m,k)$ and $\deg f_j=\deg L_j\circ F$ for $j=1,\ldots,k$, provided $\deg f_1\ge \ldots\ge \deg f_m>0$.
\end{theorem}

The above theorem gives a method of reduction of the problem of calculating the {\L}ojasiewicz exponent at infinity of overdetermined polynomial mapping to the case where the dimensions of domain and counterdomain of mapping are equal. 

\section{Proofs of Theorem \ref{thmsemialglocsets}, 
 \ref{thmsemialglocsetsinfty} and Corollary \ref{thmsemialg2}}\label{sectsemialgebr}\label{sectsemialgebrlocal}

It suffices to consider the case when $X$ and $Y$ are the basic closed semialgebraic sets. So, let
$$
X=\{x\in\rr^N:g_{1,1}(x)\ge 0,\ldots,g_{1,r(X)}(x)\ge 0,\; h_{1,1}(x)=\cdots=h_{1,l}(x)=0\}, 
$$
$$
Y=\{x\in\rr^N:g_{2,1}(x)\ge 0,\ldots,g_{2,r(Y)}(x)\ge 0,\; h_{2,1}(x)=\cdots=h_{2,l}(x)=0\}, 
$$
where $g_{i,j}, h_{i,s}\in\rr[x_1,\ldots,x_N]$. We may assume that the number of equations defining $X$ and $Y$ are equal, because we can repeat the same equations if necessary. Let $r_1=r(X)$, $r_2=r(Y)$, $r=r_1+r_2$, and let 
$G_i:\rr^N\times \rr^{r}\to\rr^{r_i}$, $i=1,2$, be the polynomial mappings defined by
$$
G_1(x,y_1,\ldots,y_r)= (g_{1,1}(x)-y_1^2,\ldots,g_{1,r_1}(x)-y_{r_1}^2),
$$
$$
G_2(x,y_1,\ldots,y_r)= (g_{2,1}(x)-y_{r_1+1}^2,\ldots,g_{2,r_2}(x)-y_{r_1+r_2}^2).
$$
 Let 
$$
A=\{(x,y_1,\ldots,y_{r})\in\rr^N\times \rr^{r}:G_{1}(x,y)=0, \; h_{1,1}(x)=\cdots=h_{1,l}(x)=0\}, 
$$
$$
B=\{(x,y_1,\ldots,y_{r})\in\rr^N\times \rr^{r}:G_{2}(x,y)=0, \; h_{2,1}(x)=\cdots=h_{2,l}(x)=0\}.
$$
Then the sets $A$ an  $B$ are algebraic and  $\pi(A)=X$, $\pi(B)=Y$, where $\pi:\rr^N\times \rr^r:(x,y)\mapsto x\in\rr^N$. Moreover,  $\deg G_1\le d$, $\deg G_2\le d$, provided $d>1$. 

By the definitions of $A$ and $B$, we immediately obtain that 
\begin{equation}\label{commonpointy}
\forall_{x_1\in X}\;\forall_{x_2\in Y}\;\exists_{y\in \rr^r}\; (x_1,y)\in A\;\land\; (x_2,y)\in B.
\end{equation}
From the definitions of the sets $A$ and $B$, we obtain that
\begin{equation}\label{eqdistance1}
\begin{split}
\forall_{x\in\rr^N\setminus X} \;\exists_{x_1\in X} \;\forall_{y\in\rr^{r}}&
\;[\dist(x,X)=|x-x_1|\,\land\,(x_1,y)\in A\;\Rightarrow \\&\; \dist(x,X) \ge \dist ((x,y),A)] 
\end{split}
\end{equation}
and
\begin{equation}\label{eqdistance2}
\begin{split}
\forall_{x\in\rr^N\setminus Y} \;\exists_{x_2\in Y} \;\forall_{y\in\rr^{r}}&
\;[\dist(x,Y)=|x-x_2|\,\land\,(x_2,y)\in B\;\Rightarrow \\&\; \dist(x,Y)\ge\dist ((x,y),B)].
\end{split}
\end{equation}
Indeed, we will prove \eqref{eqdistance1}. The proof of \eqref{eqdistance2} is similar. Take $x\in \rr^N\setminus X$ and let $x_1\in X$ satisfy $\dist (x,X)=|x-x_1|$. So, for any $y\in \rr^{r}$ 
 such that $(x_1,y)\in A$, we have
$$
\dist(x,X)=|x-x_1|=|(x,y)-(x_1,y)|\ge\dist((x,y),A).
$$
This gives \eqref{eqdistance1}.

\medskip
{\bf Proof of Theorem \ref{thmsemialglocsets}.} We will prove the assertion for non-isolated intersection $X\cap Y$ at the origin. If zero is an isolated point of the intersection $X\cap Y$, we proceed in the same way using the formula \eqref{eq71} instead of \eqref{eq7}. Let $p=d(6d-3)^{N+r-1}$.

\emph{Claim 1.} The assertion \eqref{eqsemialgsets1local} 
 is equivalent to 
\begin{equation}\label{eqlojreducedloc}
\dist (x,Y)\ge C' \dist (x,X\cap Y)^p\quad\hbox{for}\quad x \in (\partial X)\cap U_1
\end{equation}
for a neighbourhood $U_1=\{x\in \rr^N:|x|<\rho \}$ of the origin, $\rho<1$, and some positive constant $C'$, where $\partial X$ denotes the boundary of $X$ (cf. \cite[Lemma 4.2]{Cy} and \cite{KS}, proof of Theorem 2). Indeed, the implication \eqref{eqsemialgsets1local} $\Rightarrow$ \eqref{eqlojreducedloc} is obvious. Assume that the  implication  \eqref{eqsemialgsets1local}~$\Leftarrow$~\eqref{eqlojreducedloc} fails. Then for a neighbourhood $U_2=\{x\in\rr^N:|x|<\frac{\rho}{2}\}$ of the origin, there exists a sequence $a_\nu \in U_2$ such that $a_\nu\to 0$ and
\begin{equation}\label{ineqcontradictionloc}
\dist (a_\nu,X)+\dist (a_\nu,Y)< \frac{1}{\nu} \dist (a_\nu,X\cap Y)^p \quad\hbox{for}\quad \nu\in\nn.
\end{equation}
Choosing a  subsequence,
if necessary, it suffices to consider two cases: $a_\nu\not\in X$ for $\nu\in\nn$ or $a_\nu\in\Int X$ for $\nu\in\nn$.

Consider the case when $a_\nu\not\in X$ for $\nu\in\nn$. Let $x_\nu\in (\partial X)\cap U_1$ be such that $\dist (a_\nu,X)=|a_\nu-x_\nu|$. Since $\rho<1$, then we have $\dist(a_\nu,X)^{\frac{1}{p}}\ge \dist (a_\nu,X)$. So, for some $C''>0$,
$$
[\dist (a_\nu,X)+\dist (a_\nu,Y)]^{\frac{1}{p}}\ge \dist(a_\nu,X)^{\frac{1}{p}}\ge C''\dist (a_\nu,X),
$$
and, by \eqref{eqlojreducedloc}, 
$$
[\dist (a_\nu,X)+\dist (a_\nu,Y)]^{\frac{1}{p}}\ge \dist(x_\nu,Y)^{\frac{1}{p}}\ge C'' \dist (x_\nu,X\cap Y),
$$
where the above inequality is trivial if $x_\nu \in Y$. 
Since $\dist (a_\nu,X) + \dist (x_\nu,X\cap Y)\ge \dist (a_\nu,X\cap Y)$, by adding the above inequalities, we obtain 
$$
[\dist (a_\nu,X)+\dist (a_\nu,Y)]^{\frac{1}{p}}\ge \frac{C''}{2}\dist (a_\nu,X\cap Y).
$$
This contradicts \eqref{ineqcontradictionloc} and proves the Claim in the considered case.

Consider now the case when all  $a_\nu\in \Int X$.
Let $y_\nu\in Y\cap U_1$ be such that $\dist (a_\nu,Y)=|a_\nu-y_\nu|$. Then there exist $x_\nu\in(\partial X)\cap [a_\nu,y_\nu]$, where $[a_\nu,y_\nu]$ is the segment  with endpoints $a_\nu,y_\nu$.

 
 By \eqref{ineqcontradictionloc} and the choice of $\rho$, 
$$
|a_\nu-x_\nu|\le \dist (a_\nu,Y)< \frac{1}{\nu} \dist (a_\nu,X\cap Y)^p <\frac{1}{2}\dist (a_\nu,X\cap Y)\quad\hbox{for }\nu \ge 2.
$$
Hence, 
$$
\dist (x_\nu,X\cap Y)\ge \dist(a_\nu,X\cap Y)-|a_\nu-x_\nu|\ge \frac{1}{2}\dist (a_\nu,X\cap Y)\quad\hbox{for }\nu\ge 2.
$$
This and \eqref{ineqcontradictionloc} gives
$$
\dist (x_\nu,Y)\le \dist (a_\nu,Y)< \frac{1}{\nu} \dist (a_\nu,X\cap Y)^p \le \frac{2^p}{\nu} \dist (x_\nu,X\cap Y)^p \quad\hbox{for } \nu\ge 2.
$$
This contradicts \eqref{eqlojreducedloc} and proves the Claim in the considered case. Summing up we have proved Claim 1.

If $d=1$, then the assertion is trivial. Assume that $d>1$.
By \eqref{eq7}, there exists a positive constant $C$ such that
\begin{equation}\label{eqproof141}
\dist((x,y),A)+\dist((x,y),B)\ge  C\dist ((x,y),A\cap B)^{d(6d-3)^{N+r-1}} 
\end{equation}
in a neighbourhood $W$ of $0\in \rr^{N+r}$. Obviously, for any $(x,y)\in \rr^{N+r}$,
\begin{equation}\label{eqproof142}
\dist((x,y),A\cap B)\ge \dist(x,X\cap Y).
\end{equation}

One can assume that $g_{i,j}(0)=0$ for any $i,\,j$. Indeed, if $g_{i,j}(0)<0$ for some $i,j$, then $0\not\in X$ or $0\not\in Y$, which contradicts the assumption. If $g_{i,j}(0)>0$ for some $i,j$, then we can omit this inequality in the definition of $X$, respectively $Y$ and the germ at $0$ of $X$, respectrively $Y$ will not change. If $g_{i,j}(0)>0$ for any $i,j$, then the assertion reduces to \eqref{KurdykaSpodz}. So, there exists a neighbourhood $W_1=U_3\times U'\times U''\subset W$ of $0\in\rr^{N+r}$, where $U_3\subset \rr^N$, $U'\subset \rr^{r(X)}$ and $U''\subset\rr^{r(Y)}$ such that:
\begin{equation}\label{eqpoint1}
\begin{split}
\hbox{for any $(x_1,y',y'')\in A$},\;\;&\hbox{where $x_1\in \rr^N$, $y'\in \rr^{r(X)}$, $y''\in \rr^{r(Y)}$}\\
&\hbox{ if $x_1\in X\cap U_3$, then $y'\in U'$} 
\end{split}
\end{equation}
and
\begin{equation}\label{eqpoint2}
\begin{split}
\hbox{for any $(x_2,y',y'')\in B$},\;\;&\hbox{where $x_2\in \rr^N$, $y'\in \rr^{r(X)}$, $y''\in \rr^{r(Y)}$}\\
&\hbox{ if $x_2\in Y\cap U_3$, then $y''\in U''$}. 
\end{split}
\end{equation}

Let $U\subset U_3$ be a neighbourhood of $0\in\rr^N$. Take $x\in (\partial X)\cap U$, 
and let 
$x'\in Y$ be a point for which 
$\dist (x,Y)=|x-x'|$. By \eqref{commonpointy} there exists $y\in \rr^r$ such that $(x,y)\in A$ and $(x',y)\in B$. Diminishing the neighbourhood $U$, if necessary, we may assume that $x'\in U_3$. By \eqref{eqpoint1} and \eqref{eqpoint2} we see that $(x,y)\in W$, so, by \eqref{eqdistance1} and \eqref{eqdistance2}, 
$$
\dist(x,Y)\ge \dist((x,y),A)+\dist((x,y),B).
$$
Summing up, \eqref{eqproof141}, \eqref{eqproof142} and Claim 1 gives the assertion.\hfill$\square$

\medskip
{\bf Proof of Theorem \ref{thmsemialglocsetsinfty}.} Let $p=d(6d-3)^{N+r-1}$. If $X\setminus Y=\emptyset$ or $Y\setminus X=\emptyset$, then the assertion is obvious. So, we will assume that $X\setminus Y\ne \emptyset$.

By \eqref{eqKSmain} we have
\begin{equation}\label{eqKSdwa}
\dist((x,y),A)+\dist((x,y),B)\ge  C\left(\frac{\dist ((x,y),A\cap B)}{1+|(x,y)|^2}\right)^{p} 
\end{equation}
for $(x,y)\in \rr^{N+r}$. Since $\dist((x,y),A\cap B)\ge \dist(x,X\cap Y)$ for any $(x,y)\in\rr^{N+r}$ (see \eqref{eqproof142}), then \eqref{eqKSdwa} gives
\begin{equation}\label{eqKSdwa2}
\dist((x,y),A)+\dist((x,y),B)\ge  C\left(\frac{\dist (x,X\cap Y)}{1+|(x,y)|^2}\right)^{p} 
\end{equation}
for $(x,y)\in \rr^{N+r}$.

\emph{Claim 2.} The assertion \eqref{eqsemialgsets1infty} is equivalent to 
\begin{equation}\label{eqlojreduced}
\dist (x,Y)\ge C'\left(\frac{\dist (x,X\cap Y)}{1+|x|^d}\right)^p\quad\hbox{for}\quad x \in \partial X
\end{equation}
for some positive constant $C'$ 
(cf. \cite[Lemma 4.2]{Cy} and \cite{KS}, proof of Theorem 2). Indeed, the implication \eqref{eqsemialgsets1infty} $\Rightarrow$ \eqref{eqlojreduced} is obvious. Assume that the  implication  \eqref{eqsemialgsets1infty}~$\Leftarrow$~\eqref{eqlojreduced} fails. Then there exists a sequence $a_\nu \in \rr^N$ 
such that
\begin{equation}\label{ineqcontradiction}
\dist (a_\nu,X)+\dist (a_\nu,Y)< \frac{1}{\nu} \left(\frac{\dist (a_\nu,X\cap Y)}{1+|a_\nu|^d}\right)^p \quad\hbox{for}\quad \nu\in\nn.
\end{equation}
By using Theorem \ref{thmsemialglocsets} we see that $|a_\nu|\to \infty$. Choosing subsequences of the sequence $a_\nu$ if necessary, it suffices to consider two cases: $a_\nu\not\in X$ for $\nu\in\nn$ or $a_\nu\in\Int X$ for $\nu\in\nn$.

Consider the case $a_\nu\not\in X$ for $\nu\in \nn$. Let $b_\nu\in \partial X$ be such that $\dist (a_\nu,X)=|a_\nu-b_\nu|$.  Since $\left(\frac{\dist (a_\nu,X\cap Y)}{1+|a_\nu|^d}\right)^p $ is a bounded sequence then $|b_\nu-a_\nu|=\dist(a_\nu,X)\to 0$. So, for some $C''>0$ and sufficiently large $\nu$,
$$
[\dist (a_\nu,X)+\dist (a_\nu,Y)]^{\frac{1}{p}}\ge \dist(a_\nu,X)^{\frac{1}{p}}\ge C'' \left(\frac{ \dist (a_\nu,X) }{1+|a_\nu|^d}\right),
$$
and, by \eqref{eqlojreduced}, the fact that $|a_\nu|\to \infty$ and $|b_\nu-a_\nu|\to 0$,
$$
[\dist (a_\nu,X)+\dist (a_\nu,Y)]^{\frac{1}{p}}\ge \dist(b_\nu,Y)^{\frac{1}{p}}\ge C'' \left(\frac{ \dist (b_\nu,X\cap Y) }{1+|a_\nu|^d}\right),
$$
where the above inequality is trivial if $b_\nu \in Y$. 
Since $\dist (a_\nu,X) + \dist (b_\nu,X\cap Y)\ge \dist (a_\nu,X\cap Y)$, by adding the above inequalities, we obtain 
$$
[\dist (a_\nu,X)+\dist (a_\nu,Y)]^{\frac{1}{p}}\ge\frac{C''}{2}\left(\frac{ \dist (a_\nu,X\cap Y) }{1+|a_\nu|^d}\right).
$$
This contradicts \eqref{ineqcontradiction} and proves the Claim in the considered case.

Consider the case $a_\nu \in \Int X$ for $\nu\in\nn$. Let $y_\nu\in Y$ be such that $\dist (a_\nu,Y)=|a_\nu-y_\nu|$. Then there exist $x_\nu\in(\partial X)\cap [a_\nu,y_\nu]$ for $\nu\in\nn$. By \eqref{ineqcontradiction}, for sufficiently large $\nu$
\begin{equation}\label{eqClaim2}
|a_\nu-x_\nu|\le \dist (a_\nu,Y)< \frac{1}{\nu} \left(\frac{\dist (a_\nu,X\cap Y)}{1+|a_\nu|^d}\right)^p  <\frac{1}{2} \dist (a_\nu,X\cap Y).
\end{equation}
Hence, 
$$
\dist (x_\nu,X\cap Y)\ge \dist(a_\nu,X\cap Y)-|a_\nu-x_\nu|\ge \frac{1}{2}\dist (a_\nu,X\cap Y).
$$
This and \eqref{ineqcontradiction} gives
$$
\dist (x_\nu,Y)\le \dist (a_\nu,Y)< \frac{1}{\nu} \left(\frac{\dist (a_\nu,X\cap Y)}{1+|a_\nu|^d}\right)^p \le \frac{2^p}{\nu} \left(\frac{\dist (x_\nu,X\cap Y)}{1+|a_\nu|^d}\right)^p.
$$
By \eqref{eqClaim2}, for sufficiently large $\nu$, $|x_\nu|\le 2|a_\nu|$, so, for a positive constant $C'''$,
$$
\dist (x_\nu,Y)\le \frac{C'''}{\nu} \left(\frac{\dist (x_\nu,X\cap Y)}{1+|x_\nu|^d}\right)^p.
$$
This contradicts \eqref{eqlojreduced} and proves the Claim in the considered case. Summing up we have proved Claim 2.

Take any $x_0\in {\partial X}$. By \eqref{commonpointy}, \eqref{eqdistance1} and \eqref{eqdistance2} there exist $x_2\in Y$ and $y_0\in\rr^r$ such that $(x_0,y_0)\in A$, $(x_2,y_0)\in B$, and 
$\dist (x_0,Y)=|x_0-x_2|\ge \dist ((x_0,y_0),B)$. Hence from \eqref{eqKSdwa2},
\begin{equation}\label{eqKSdwa3}
\dist(x_0,Y)\ge  C\left(\frac{\dist (x_0,X\cap Y)}{1+|(x_0,y_0)|^2}\right)^{p} .
\end{equation}

It is easy to observe that there exist constants $C_1,R_1>0$ such that for $(x,y)\in A$, $|(x,y)|\ge R_1$ we have $C_1|y|^2\le |x|^d$. 
 Since $d\ge 2$, then for a constant $C_2>0$ we obtain $|(x,y)|\le C_2|x|^{d/2}$ for $(x,y)\in A$, $|(x,y)|\ge R_1$. Hence from \eqref{eqKSdwa3} we easily deduce  
\begin{equation}\label{eqKSdwa4}
\dist(x_0,Y)\ge  C\left(\frac{\dist (x_0,X\cap Y)}{1+C_2^2|x_0|^d}\right)^{p},
\end{equation}
provided $|x_0|\ge R_1$. So, diminishing $C$, if necessary, we obtain  
 \eqref{eqKSdwa4} for arbitrary $x_0\in \partial{X}$. 
This, together with the Claim 2 gives the assertion of Theorem \ref{thmsemialg}.\hfill$\square$

\medskip
{\bf Proof of Corollary \ref{thmsemialg2}}
Let $H:\rr^{N+r}\to \rr^{m+r+l}$ be a polynomial mapping defined by
$$
H(x,y)=(F(x),G_1(x,y), h_{1,1}(x), \ldots,h_{1,l}(x)),\quad x\in\rr^N, y\in\rr^r.
$$
Then $\deg H\le d$. Let $V=F^{-1}(0)\cap X$ and let $Z=H^{-1}(0)$. 
By \eqref{eq8KS}, for some positive constant $C$, we have
$$
|H(x,y)|\ge C\left(\frac{\dist((x,y),Z)}{1+|(x,y)|^2}\right)^{{d(6d-3)^{N+r-1}}} \quad\hbox{for}\quad (x,y)\in\er^N\times\er^r.
$$
Because $\dist((x,y),Z)\ge \dist (x,V)$, then
\begin{equation}\label{eqproofKS1}
|H(x,y)|\ge C\left(\frac{\dist(x,V)}{1+|(x,y)|^2}\right)^{{d(6d-3)^{N+r-1}}} \quad\hbox{for}\quad (x,y)\in\er^N\times\er^r.
\end{equation}
It is easy to observe that there exist constants $C_1,R_1>0$ such that for $(x,y)\in A$, $|(x,y)|\ge R_1$ we have $C_1|y|^2\le |x|^D$. 
 Since $D\ge 2$, then for a constant $C_2>0$ we obtain $|(x,y)|\le C_2|x|^{D/2}$ for $(x,y)\in A$, $|(x,y)|\ge R_1$. Hence from \eqref{eqproofKS1} we easily deduce  
\eqref{eqsemialg1KS2} for $x\in X$, $|x|\ge R_1$. So, diminishing $C$, if necessary, we obtain  
 \eqref{eqsemialg1KS2} for arbitrary $x\in X$

We will show the second part of the assertion. Since $X$ is an unbounded set, we may assume that the set $A$ is unbounded, too. 
Since $V$ is a compact set, then the set $H^{-1}(0)$ is also compact. By \eqref{eq83} we have $\wykl_\infty (H)\ge -d(6d-3)^{N+r-1}$, in particular for some constants $C,R>0$,
\begin{equation}\label{eqproofsemialg}
|H(x,y)|\ge C|(x,y)|^{-d(6d-3)^{N+r-1}}\quad \hbox{for}\quad (x,y)\in A,\quad |(x,y)|\ge R.
\end{equation}
Since $|(x,y)|\le C_2|x|^{D/2}$ for $(x,y)\in A$, $|(x,y)|\ge R_1$, so, for some constant $C_3>0$. 
$$
|F(x)|=|H(x,y)|\ge C_3|x|^{-\frac{D}{2}d(6d-3)^{N+r-1}} \quad\hbox{for}\quad (x,y)\in A,\quad |(x,y)|\ge R,
$$
and $\wykl^\rr_\infty(F|X)\ge -\frac{D}{2}d(6d-3)^{N+r-1}$. 
This ends the proof of Corollary \ref{thmsemialg2}. \hfill$\square$

\section{Proof of Theorem \ref{Theorem1}}\label{Proofssection}


 Let $k\in\ez$, $n\le k\le m$. Take a closed semialgebraic set $Z\subset \rr^N$ of dimension $\dim_{\rr}Z=n$, and let 
 $$
 \pi :Z\ni(x,y)\mapsto y\in\rr^m.
 $$
 Then the set $\pi(Z)$ is semialgebraic with $\dim_\rr \pi(Z)\le n$. Denote by $Y\subset \cc^m$ the complex Zariski closure of $\pi(Z)$. So, $Y$ is an algebraic set of complex dimension $\dim_{\cc}Y\le n$.

Assume that $0\in Y$. Let $C_0(Y)\subset \cc^m$ be the tangent cone to $Y$ at $0$ in the sense of Whitney \cite[p. 510]{w}. It is known that $C_0(Y)$ is an algebraic set and $\wym_\ec C_0(Y)\leq n$. So, we have

\begin{lemma}\label{Lemma 1.}
 For the generic $L\in \LL^\kk(m,k)$, 
$$
L^{-1}(0)\cap C_0(Y)\subset \{0\}.
$$
\end{lemma}

In the proofs of Theorems \ref{Theorem1}, \ref{Theorem2} and \ref{Theorem3} we will need the following

\begin{lemma}\label{Lemma 2.}
If $L\in \LL^\kk(m,k)$ satisfies $L^{-1}(0)\cap C_0(Y)\subset \{0\}$, then there exist $\varepsilon,\,C_1,\,C_2>0$ such that for $z\in Z$, $|\pi(z)| <\varepsilon$ we have 
\begin{equation}\label{eqlemma2}
C_1 |\pi(z)| \leq | L(\pi(z))| \leq C_2 |\pi(z)|.
\end{equation}
\end{lemma}

{\bf Proof.} It is obvious that for $C_2=\parallell L\parallelp$ we obtain $ | L(\pi(z))| \leq C_2|\pi(z)|$ for $z\in Z$. This gives the right hand side inequality in \eqref{eqlemma2}.

Now, we show the left hand side inequality in \eqref{eqlemma2}. Assume to the contrary, that for any $\varepsilon,\,C_1>0$ there exists $z\in Z$ such that 
$$
C_1 | \pi(z) | > | L(\pi(z)) |\quad \textrm{ and }\quad | \pi(z) | <\varepsilon .
$$ 
In particular, for $\nu \in \en$, $C_1=\frac{1}{\nu }$, $\varepsilon =\frac{1}{\nu }$ there exists $z_\nu \in Z$ such that 
$$
\frac{1}{\nu }|\pi(z_\nu )| > | L(\pi(z_\nu ))| \quad\textrm{ and }\quad | \pi(z_\nu) | <\frac{1}{\nu} .
$$
Thus $|\pi(z_\nu)|>0$ and
\begin{equation}\label{w1}
\frac{1}{\nu} > \frac{1}{|\pi(z_\nu )|}| L(\pi(z_\nu ))| = \left|L\left(\frac{1}{|\pi(z_\nu )|} \pi(z_\nu )\right)\right|.
\end{equation}
Let $\lambda _\nu =\frac{1}{\mid \pi(z_\nu )\mid }$ for $\nu \in \en$. Then $| \lambda_\nu \pi(z_\nu ) | = 1$
so, by choosing subsequence, if necessary, we may assume that $\lambda _\nu \pi(z_\nu )\rightarrow v$ when $\nu \rightarrow \infty $, where $v\in \cc^m$, $| v| =1$ and $\pi(z_\nu )\rightarrow 0$ as $\nu \rightarrow \infty $, thus $v\in C_0(Y)$ and $v \not = 0$. Moreover, by (\ref{w1}), we have $L(v)=0$.  So $v\in L^{-1}(0)\cap C_0(Y)\subset \{0\}.$ This contradicts the assumption and ends the proof.\hfill$\square$

\medskip
We will also need the following lemma (cf. \cite{P2}, \cite{Sp3}). Let $X\subset \rr^N$ be a closed semialgebraic set such that $0\in X$.

\begin{lemma}\label{Lemma3} Let $F,\,G:(\er^N,0)\to (\er^m,0)$ be analytic mappings, such that\linebreak $\ord_0(F-G)>\wykl_0^\rr(F|X)$. If $0$ is an isolated zero of $F|X$ then $0$ is an isolated zero of $G|X$ and for some positive constants $\varepsilon,C_1,C_2$,
\begin{equation}\label{eq3prooftheorem1}
C_1|F(x)|\le |G(x)|\le C_2|F(x)|\quad\hbox{for}\quad x\in X,\quad |x|<\varepsilon.
\end{equation}
In particular, $\wykl_0^\rr(F|X)=\wykl_0^\rr(G|X)$.
\end{lemma}

{\bf Proof.} Since $F$ is a Lipschitz mapping, so $1\le \wykl_0^\rr(F|X)<\infty$ and for some positive constants $\varepsilon_0,C$,
\begin{equation}\label{eq1prooflem3}
|F(x)|\ge C|x|^{\wykl_0^\rr(F|X)}\quad \hbox{for}\quad x\in X,\quad |x|<\varepsilon_0.
\end{equation}
By the assumption $\ord_0 (F-G)>\wykl_0^\rr(F|X)$ it follows that there exist $\eta\in\er$,\linebreak $\eta>\wykl_0^\rr(F|X)$ and $\varepsilon_1>0$ such that  $||F(x)|-|G(x)||\le |x|^\eta$ for $x\in X$, $|x|<\varepsilon_1$. Assume that \eqref{eq3prooftheorem1} fails. Then for some sequence $x_\nu\in X$ such that $x_\nu \to 0$ as $\nu\to\infty$, we have 
$$
\frac{1}{\nu}|F(x_\nu)|>|G(x_\nu)|\quad \hbox{or}\quad \frac{1}{\nu}|G(x_\nu)|> |F(x_\nu)|\quad\hbox{for}\quad \nu\in\en.
$$
So, in the both above cases, by \eqref{eq1prooflem3} for $\nu\ge 2$, we have
$$
\frac{C}{2}|x_\nu|^{\wykl_0^\rr(F|X)}\le \frac{1}{2}|F(x_\nu)|< |F(x_\nu)-G(x_\nu)|\le |x_\nu|^\eta,
$$
which is impossible. The last part of the assertion follows immediately from \eqref{eq3prooftheorem1}.

\hfill$\square$

\medskip
{\bf Proof of Theorem \ref{Theorem1}.} The assertion \eqref{eq1realisolated} we prove analogously as Theorem 2.1 in \cite{Sp3}. We will prove the second part of the assertion.

Let $G=(g_1,\ldots,g_m):(\er^N,0)\to(\er^m,0)$ be a polynomial mapping such that $\ord_0^\rr(F-G)>\wykl_0^\rr(F|X)$. Obviously, such a mapping $G$ does exist. 
By Lemma \ref{Lemma3}, $\wykl_0^\rr(F|X)=\wykl_0^\rr(G|X)$ and $0$ is an isolated zero of $G|X$. Taking, if necessary, intersection of X with a ball $B$ centered at zero, we can assume that $\dim_{\rr,0} X=\dim_\rr X$. So, by Lemmas \ref{Lemma 1.} and \ref{Lemma 2.} for the generic $L\in\LL^\rr(m,k)$ we have that $L\circ G|X$ has an isolated zero at $0\in\er^n$, $\wykl_0^\rr(G|X)=\wykl_0^\rr(L\circ G|X)$, and
\begin{equation*}
\begin{split}
\ord_0(L\circ G-L\circ F)=&\ord_0 L\circ(G-F)\ge \ord_0 (G-F)\\
>&\wykl_0^\rr (F|X)=\wykl_0^\rr(G|X)=\wykl_0^\rr(L\circ G|X),
\end{split}
\end{equation*}
so, by Lemma \ref{Lemma3}, $\wykl_0^\rr(L\circ F|X)=\wykl_0^\rr(L\circ G|X)=\wykl_0^\rr(F|X)$. This gives \eqref{eq2realisolated}. The particular part of the assertion is proved analogously as in \cite[Proposition 2.1]{Sp3}.
\hfill$\square$

\section{Proof of Theorem \ref{Theorem2}}\label{Proofssection2}

Let $X\subset \rr^N$ be a closed semialgebraic set $\dim_\rr X=n$, and let $0\in X$.  Taking, if necessary, intersection of X with a ball $B$ centered at zero, we can assume that $\dim_{\rr,0} X=\dim_\rr X$. 

%

From \cite[Proposition 1.1]{Sp2} we immediately obtain

\begin{proposition}\label{Proposition1.1}  Let $G=(g_1,...,g_m):X\to \kk^m$ be a semialgebraic 
mapping, $g_j\ne 0$ for $j=1,...,m$, where $m\ge n\ge 1$, and let $k\in\mathbb{Z}$, $n\le k\le m$.

(i) For the generic $L\in \LL^\kk(m,k)$,
\begin{equation}\label{eq1tw3}
\#[(L\circ G)^{-1}(0)\setminus G^{-1}(0)]<\infty.
\end{equation}

(ii) For the generic $L\in \Ld^\kk(m,k)$,
\begin{equation}\label{eq1tw3ii}
\#[(L\circ G)^{-1}(0)\setminus G^{-1}(0)]<\infty.
\end{equation}
%
\end{proposition}

{\bf Proof.} Let $Y\subset \cc^N\times \cc^m$ be the Zariski closure of the graph of $G$, and let $\pi:Y\ni(x,y)\mapsto y\in \cc^m$. Then for $(x,y)\in Y$ such that $x\in X$ and $y\in \kk^m$ we have $y=G(x)$. Let us consider the case $n=k$. Let 
$$
U=\{L\in \LL^\cc (m,n):\# [(L\circ \pi)^{-1}(0)\setminus \pi^{-1}(0)]<\infty \}.
$$ 
By Proposition 1.1 in \cite{Sp2}, $U$ contains a non-empty Zariski open subset of $\LL^\cc(m,n)$. Then $U$ contains a dense Zariski open subset $W$ of $\LL^\rr(m,n)$. This gives the assertion (i) in the case $n=k$. 

Let now $k>n$. Since for $L=(L_1,\ldots,L_k)\in \LL^\kk(m,k)$,
$$
(L\circ \pi)^{-1}(0)\subset ((L_1,\ldots, L_n)\circ \pi)^{-1}(0),
$$
then the assertion (i) follows from the previous case. We prove the assertion (ii) analogously as \cite[Proposition 1.1]{Sp2}.
\hfill$\square$

\medskip
{\bf Proof of Theorem \ref{Theorem2}.} 
Without loss of generality we may assume that $F\ne 0$. By the definition, there exist $C,\varepsilon>0$ such that for $x\in X$, $|x|<\varepsilon$ we have
\begin{equation}\label{eq1proof2}
|F(x)|\ge C\dist (x,F^{-1}(0))^{\wykl^\kk_0(F|X)},
\end{equation}
and $\wykl_0^\kk(F|X)$ is the smallest exponent for which the inequality holds. 
Let $L\in \LL^\kk(m,k)$ be such that $F^{-1}(0)\cap U_L=(L\circ F)^{-1}(0)\cap U_L$ for some neighbourhood $U_L\subset \kk^N$ of $0$. Diminishing $\varepsilon$ and the neighbourhood $U_L$, if necessary, we may assume that the equality $\dist(x,F^{-1}(0))=\dist(x,F^{-1}(0)\cap U_L)$ holds for $x\in X$, $|x|<\varepsilon$. Obviously $L\ne 0$, so, $||L||>0$, and $|F(x)|\ge \frac{1}{||L||}|L(F(x))|$. Then by \eqref{eq1proof2} we obtain $\wykl^\kk_0(F|X)\le\wykl^\kk_0(L\circ F|X)$, and \eqref{eq1lokalnonisolated} is proved.

By Proposition \ref{Proposition1.1} and Lemmas \ref{Lemma 1.} and \ref{Lemma 2.}, for the generic $L\in \LL^\kk(m,k)$ we have that $F^{-1}(0)\cap U_L=(L\circ F)^{-1}(0)\cap U_L$ for some neighbourhood $U_L\subset \kk^N$ of $0$ and 
there exist $\varepsilon,\,C_1,\,C_2>0$ such that for $x\in X$, $| x | <\varepsilon$, 
\begin{equation}\label{eq2proof2}
C_1 |F(x)| \leq | L(F(x))| \leq C_2 | F(x)|.
\end{equation}
This and \eqref{eq1proof2} gives \eqref{eq2lokalnonisolated} and ends the proof of Theorem \ref{Theorem2}.\hfill$\square$

\section{Proof of Theorem \ref{Theorem3}}\label{ProofofTheorem}

The argument of Lemma 2.2 from \cite{Sp2} gives

\begin{lemma}\label{lemma4}
Let $F:X\to\rr^m$ with $m\ge n=\dim_\rr X$ be a semialgebraic mapping, where $X\subset \rr^N$, and let $n\le k\le m$. Then there exists a Zariski open and dense subset $U\subset \LL^\rr(m,k)$ 
such that for any $L\in U$ and any $\varepsilon>0$ there exist $\delta>0$ and $r>0$ such that for any $x\in X$,
$$
|x|>r\;\land \;|L\circ F(x)|<\delta\;\Rightarrow\;|F(x)|<\varepsilon.
$$
\end{lemma}

{\bf Proof} (cf. proof of Lemma 2.2 in \cite{Sp2}).  Let us consider the case $k=n$. Let $W\subset \cc^N$ be the Zariski closure of $F(X)$. Then $\dim_\cc W\le n$. In the case $\dim_\cc W<n$, by Lemma 2.1 in \cite{Sp2} we easily obtain the assertion. Assume that $\dim W=n$. We easily see that for an algebraic set $V\subset W$, $\dim_\cc V\le n-1$, the mapping 
$F|_{X\setminus F^{-1}(V)}:X\setminus F^{-1}(V)\to W\setminus V$ is  proper. By Lemma 2.1 in \cite{Sp2} there exists a Zariski open and dense subset $U_1\subset \LL^\rr(m,k)$ such that for any $L\in U_1$ and for any $\varepsilon>0$ there exists $\delta>0$ such that for $z\in V$,
\begin{equation}\label{eqSp2Le4ma211}
|L(z)|<\delta \Rightarrow |z|<\varepsilon.
\end{equation}
 Moreover, for $L\in U_1$,
\begin{equation}\label{eqSp2Lema21}
W\subset \{z\in \cc^m:|z|\le C_L(1+|L(z)|)\}
\end{equation}
for some $C_L>0$. 

 Let 
 $$
 U=\{L\in \LL^\rr(m,n):L\in U_1\}.
 $$
 Obviously, $U$ is a dense and Zariski open subset of $\LL^\rr(m,n)$. 
  Take $L\in U$ and $\varepsilon >0$. Assume to the contrary that there exists a sequence $x_\nu\in X$ such that $|x_\nu|\to \infty$, $|L(f(x_\nu))|\to 0$ and $|f(x_\nu)|\ge \varepsilon$. By \eqref{eqSp2Lema21} we may assume that $f(x_\nu)\to y_0$ for some $y_0\in W$. Since $F|_{X\setminus F^{-1}(V)}:X\setminus F^{-1}(V)\to W\setminus V$ is a proper mapping, we have $y_0\in V$. So, $|y_0|\ge \varepsilon$ and $L(y_0)=0$. This contradicts \eqref{eqSp2Le4ma211} and ends the proof in the case $n=k$.
  
Let now, $k>n$ and let
$$
U=\{L=(L_1,\ldots,L_k)\in\LL^\rr(m,k):(L_1,\ldots,L_n)\in U_1\}. 
$$
Then for any $L=(L_1,\ldots,L_k)\in U$ and $x\in \er^n$ we have 
$$
|(L_1,\ldots,L_n)\circ F(x)|\le |L\circ F(x)|,
$$
so, the assertion immediately follows from the previous case. 
\hfill$\square$

\medskip

{\bf Proof of Theorem \ref{Theorem3}} (cf. proof of Theorem 2.1 in \cite{Sp2}). Since for non-zero $L\in\LL^\rr(m,k)$ we have $|L\circ F(x)|\le ||L|||F(x)|$ and $||L||>0$,  then by the definition of the {\L}ojasiewicz exponent at infinity we obtain the first part of the assertion. We will prove the second part of the assertion.

Since $F^{-1}(0)$ is a compact set, by Proposition \ref{Proposition1.1}, there exists a dense Zariski open subset $U $ of $\LL^\rr(m,k)$ such that 
$$
U\subset \{L\in\LL^\rr(m,k): (L\circ F)^{-1}(0)\hbox{ is a compact set}\}.
$$
So, for the generic $L\in\LL^\rr(m,k)$ the set $(L\circ F)^{-1}(0)$ is compact.

If $\wl^\er (F|X)< 0$, the assertion (\ref{eq2global}) follows from Lemmas \ref{Lemma 1.}, \ref{Lemma 2.} and \ref{lemma4}.

Assume that $\wl^\rr(F|X)=0$. Then there exist $C,\, R>0$ such that $|F(x)|\ge C$ as $|x|\ge R$. Moreover, there exists a sequence $x_\nu\in X$ such that $|x_\nu|\to \infty$ as $\nu\to\infty$ and $|F(x_\nu)|$ is a bounded sequence. So, by Lemma \ref{lemma4} for the generic $L\in U$ 
 and $\varepsilon =C$ there exist $r,\delta >0$ such that $|L\circ F(x)|\ge \delta$ as $|x|>r$, so $\wl^\rr(L\circ F|X)\ge 0$. Since $|L\circ F(x_\nu)|$ is a bounded sequence, we have $\wl^\rr(L\circ F|X)\le 0$. Summing up $\wl^\rr(L\circ F|X)=\wl^\rr(F|X)$ in the considered case.

In the case $\wl^\er (F|X)>0$, we obtain the assertion analogously as in the proof of Theorem 2.1 in \cite{Sp2}.\hfill$\square$


{\small
\begin{thebibliography}{99}

\normalsize

\bibitem{ATW}
R. Achilles, P. Tworzewski, T. Winiarski, \emph{\it On improper isolated intersection in complex analytic geometry},  Ann. Polon. Math.  51  (1990), 21--36.





\bibitem{BPR1}  S. Basu, R. Pollack, M-F. Roy, \emph{Computing the dimension of a semi-algebraic set}. Zap. Nauchn. Sem. S.-Peterburg. Otdel. Mat. Inst. Steklov. (POMI) 316 (2004), Teor. Slozhn. Vychisl. 9, 42--54, 225; translation in J. Math. Sci. (N. Y.) 134 (2006), no. 5, 2346–2353.

\bibitem{BPR2} S. Basu, R. Pollack, M-F. Roy, \emph{On the combinatorial and algebraic complexity of quantifier elimination}. J. ACM 43 (1996), no. 6, 1002–1045.


\bibitem{BPR} S. Basu, R. Pollack, M-F. Roy, \emph{Algorithms in real algebraic geometry. Second edition. Algorithms and Computation in Mathematics}, Vol 10. Springer-Verlag, Berlin, 2006. 662 pp.

\bibitem{BCR} J. Bochnak, M. Coste, M-F. Roy, \emph{Real algebraic geometry},  Springer-Verlag, Berlin, 1998. 



\bibitem{BR} J. Bochnak, J. J. Risler, \emph{Sur les exposants de {\L}ojasiewicz}, Comment. Math. Helv. 50 (1975), 493-507.


\bibitem{Brocker1} L. Br\"ocker, \emph{Minimale Erzeugung von Positivbereichen}. Geom.Dedicata 16 (1984), no. 3, 335–350.

\bibitem{Brocker0} L. Br\"ocker, \emph{On basic semialgebraic sets}. Exposition. Math. 9 (1991), no. 4, 289–334.



\bibitem{CK3} J. Ch\c adzy\'nski, \emph{On proper polynomial mappings}, Bull. Polish Acad. Sci. Math. 31 (1983), 115--120. 




\bibitem{Cy} E. Cygan, \emph{A note on separation of algebraic sets and the {\L}ojasiewicz exponent for polynomial mappings},  Bull. Sci. Math.  129  (2005),  no. 2, 139--147.

\bibitem{Cy2}E. Cygan, \emph{Intersection theory and separation exponent in complex analytic geometry}, Ann. Polon.
Math. 69 (3) (1998) 2870299.

\bibitem{CKT} E. Cygan, T. Krasi\'nski, P. Tworzewski, \emph{Separation of algebraic sets and the {\L}ojasiewicz exponent of polynomial mappings},  Invent. Math.  136  (1999),  no. 1, 75--87.






\bibitem{J1} Z. Jelonek, \emph{On the effective Nullstellensatz},  Invent. Math.  162  (2005),  no. 1, 1--17. 

\bibitem{J2} Z. Jelonek, \emph{On the {\L}ojasiewicz exponent}, Hokkaido Math. J.  35  (2006),  no. 2, 471--485. 


\bibitem{JKS} S. Ji, J. Koll\'ar, B. Shiffman, \emph{A global {\L}ojasiewicz
inequality for algebraic va\-rie\-ties}, Trans. Amer. Math. Soc.
 329 (1992), 813--818.

\bibitem{K1}J. Koll\'ar, \emph{Sharp effective Nullstellensatz}, 
J. Amer. Math. Soc.  1  (1988),  963--975.

\bibitem{K2} J. Koll\'ar, \emph{An effective {\L}ojasiewicz inequality for real polynomials},  Period. Math. Hungar.  38  (1999),  no. 3, 213--221. 


    
    

\bibitem{KMP} K. Kurdyka, T. Mostowski, A. Parusi\'nski, \emph{Proof of the gradient conjecture of R. Thom}, 
Ann. of Math. (2) 152 (2000), no. 3, 763--792.

\bibitem{KS} K. Kurdyka, S. Spodzieja, \emph{Separation of real algebraic sets and the {\L}ojasiewicz exponent}. Proc. Amer. Math. Soc. 142 (2014), no. 9, 3089–3102.


\bibitem{LejeuneTeissier}
M. Lejeune-Jalabert, B. Teissier, \emph{Cl\^{o}ture int\'egrale des id\'eaux et \'equisingularit\'e}. Centre de Math\'ematiques Ecole Polytechnique Palaiseau, 1974.















\bibitem{P2}A. P{\l}oski, \emph{Multiplicity and the {\L}ojasiewicz exponent},  
  Banach Center Publications 20, Warsaw (1988),  353 -- 364.



\bibitem{RS1}
T. Rodak, S. Spodzieja, \emph{Effective formulas for the {\L}ojasiewicz exponent at infinity},  J. Pure Appl. Algebra  213  (2009), 1816--1822.

\bibitem{RS2}
T. Rodak, S. Spodzieja, \emph{Effective formulas for the local {\L}ojasiewicz exponent}, Math. Z. 268 (2011), 37-44.

\bibitem{RS3}
T. Rodak, S. Spodzieja, \emph{{\L}ojasiewicz exponent near the fibre of a mapping}, Proc. Amer. Math. Soc. 139 (2011), 1201-1213. 

\bibitem{RS4} T. Rodak, S. Spodzieja,  \emph{Equivalence of mappings at infinity}. Bull. Sci. Math. 136 (2012), no. 6, 679–686. 

\bibitem{Roy} M-F. Roy, N. Vorobjov, \emph{The complexification and degree of a semi-algebraic set}. Math. Z. 239 (2002), no. 1, 131–142. 







\bibitem{Scheiderer} C. Scheiderer, \emph{Stability index of real varieties}. Invent. Math. 97 (1989), no. 3, 467–483.

\bibitem{Sp2} S. Spodzieja, \emph{The {\L}ojasiewicz exponent at infinity for overdetermined 
polynomial mappings},  Ann. Polon. Math.  78  (2002),  1--10. 

\bibitem{Sp3} S. Spodzieja, \emph{Multiplicity and the {\L}ojasiewicz
exponent}, Ann. Polon. Math. 73 (2000), 257--267

\bibitem{Sp4} S. Spodzieja, \emph{The Łojasiewicz exponent of subanalytic sets}, Ann. Polon. Math. 87 (2005), 247–263.

\bibitem{SS} 
S. Spodzieja, A. Szlachcińska \emph{{\L}ojasiewicz exponent of overdetermined mappings}. Bull. Pol. Acad. Sci. Math. 61 (2013), no. 1, 27–34.

\bibitem{T} B. Teissier, \emph{Vari\'et\'es polaires. I. Invariants polaires des singularit\'es d'hypersurfaces}, Invent. Math. 40 (1977), no. 3, 267--292.

\bibitem{Tw} P. Tworzewski, \emph{Intersection theory in complex analytic geometry}, Ann. Polon. Math. 62 (1995), 177-191. 

\bibitem{w} H. Whitney, \textit{Tangents to an analytic variety}, Ann. of Math. 81 (1965), 496-549.

\end{thebibliography}
}
\end{document}